\documentclass[12pt]{amsart}

\usepackage{amsfonts, amsthm, amsmath}
\usepackage[dvips]{graphics}

\theoremstyle{plain}
\newtheorem{theorem}{Theorem}
\newtheorem{question}[theorem]{Question}
\newtheorem{lemma}[theorem]{Lemma}
\newtheorem{corollary}[theorem]{Corollary}

\begin{document}

\title{Embedded Minimal Disks with Prescribed Curvature Blowup}
\author{Brian Dean}
\address{Department of Mathematics\\
    Hylan Building\\
    University of Rochester\\
    Rochester, NY  14627}
\email{bdean@math.rochester.edu}
\date{}
\thanks{The author thanks W. Minicozzi for his many helpful discussions.}

\begin{abstract}We construct a sequence of compact embedded
minimal disks in a ball in $\mathbb{R}^3$, whose boundaries lie in
the boundary of the ball, such that the curvature blows up only at
a prescribed discrete (and hence, finite) set of points on the
$x_3-\mbox{axis}$. This extends a result of Colding and Minicozzi,
who constructed a sequence for which the curvature blows up only
at the center of the ball, and is a partial affirmative answer to
the larger question of the existence of a sequence for which the
curvature blows up precisely on a prescribed closed set on the
$x_3-\mbox{axis}$.
\end{abstract}

\maketitle

In \cite{CM1}, T.H. Colding and W.P. Minicozzi II constructed a
sequence of compact embedded minimal disks in a ball in
$\mathbb{R}^3$, with boundaries lying in the boundary of the ball,
such that the curvature blows up only at the center.  This result
raises the following question.

\begin{question}\label{closedsetquestion}
Does there exist a sequence of compact embedded minimal disks in a
ball in $\mathbb{R}^3$, whose boundaries lie in the boundary of
the ball, such that the curvature blows up precisely on a
prescribed closed set on the $x_3-\mbox{axis}$?
\end{question}

Beyond that, it is interesting to consider which curves can arise as the singular set for curvature of a sequence of embedded minimal disks.  W. Meeks and M. Weber have constructed examples (see \cite{MW}) in which the singular set is a circle.

By scaling, it suffices to consider embedded minimal disks in the
unit ball.  Our main result says that the answer to Question~\ref{closedsetquestion}
is affirmative in the case where the closed set is a discrete (and
hence, finite) set of points.

\begin{theorem}\label{finitepointstheorem}
Given $n$ points $(0,0,b_j)\subset B_1,
\,b_1<\ldots <b_n$, there is a sequence of compact embedded
minimal disks $0\in\Sigma_i\subset B_1\subset\mathbb{R}^3$ with
$\partial\Sigma_i\subset\partial B_1$ and containing the vertical
segment $\{(0,0,t): |t|<1\}\subset\Sigma_i$, and such that the
following hold:
\begin{enumerate}
  \item[(i)]
  $\lim_{i\rightarrow\infty}|A_{\Sigma_i}|^2(0,0,b_j)=\infty,\,j=1,\ldots,n$
  \item[(ii)]
  $\sup_i\sup_{\Sigma_i\verb"\"\cup_jB_{\delta_j}(0,0,b_j)}|A_{\Sigma_i}|^2<\infty$
  for all $\delta_j>0,\,j=1,\ldots,n$
  \item[(iii)]
  $\Sigma_i\verb"\"\{x_3-\mbox{axis}\}=\Sigma_{1,i}\cup\Sigma_{2,i}$
  for multi-valued graphs $\Sigma_{1,i}$ and $\Sigma_{2,i}$
  \item[(iv)] $\Sigma_i\verb"\"\cup_j\{x_3=b_j\}$ converges to $n+1$
  embedded minimal disks $\Sigma^k, \,k=1,\ldots,n+1$,
  satisfying the following:
  \begin{enumerate}
    \item[(a)] $\Sigma^k\subset\{x_3<b_k\}$ for $k=1,\ldots,n$,
    and $\Sigma^{n+1}\subset\{x_3>b_n\}$
    \item[(b)] $\overline{\Sigma^1}\verb"\"\Sigma^1
    =B_1\cap\{x_3=b_1\},\,\overline{\Sigma^{n+1}}\verb"\"\Sigma^{n+1}
    =B_1\cap\{x_3=b_n\}$, and for
    $k=2,\ldots,n$,\,$\overline{\Sigma^k}\verb"\"\Sigma^k
    =B_1\cap(\{x_3=b_{k-1}\}\cup\{x_3=b_k\})$
    \item[(c)]
    $\Sigma^1\verb"\"\{x_3-\mbox{axis}\}=\Sigma_1^1\cup\Sigma_2^1$
    for multi-valued graphs $\Sigma_1^1$ and $\Sigma_2^1$ each of
    which spirals into $\{x_3=b_1\}$.
    $\Sigma^{n+1}\verb"\"\{x_3-\mbox{axis}\}=\Sigma_1^{n+1}\cup\Sigma_2^{n+1}$
    for multi-valued graphs $\Sigma_1^{n+1}$ and $\Sigma_2^{n+1}$
    each of which spirals into $\{x_3=b_n\}$.  For $k=2,\ldots,n,
    \,\Sigma^k\verb"\"\{x_3-\mbox{axis}\}=\Sigma_1^k\cup\Sigma_2^k$
    for multi-valued graphs $\Sigma_1^k$ and $\Sigma_2^k$ each of
    which spirals into $\{x_3=b_{k-1}\}$ and $\{x_3=b_k\}$.
  \end{enumerate}
\end{enumerate}
\end{theorem}

Question~\ref{closedsetquestion} remains open for closed sets in
general; for example, closed intervals or Cantor-type sets.  We
conjecture that the answer is affirmative in general.  Let us briefly discuss how one might show this.  One idea is to note that, given any closed set, there exists a countable dense subset.  One would want to construct a sequence of compact embedded minimal disks whose curvature blows up on the countable dense subset.  By Theorem~\ref{finitepointstheorem}, such a sequence exists for any set of $n$ points, for any fixed finite $n$; one would then want to let $n\rightarrow\infty$ and use a diagonal argument to obtain a sequence whose curvature blows up on the countable set of points in the dense subset.  As a result of \cite[Lemma I.1.4]{CM2}, the set of points on which the curvature blows up must be closed.  Hence, the curvature would blow up on the closure of the countable dense subset, which is precisely our prescribed closed set.

The key to extending Theorem~\ref{finitepointstheorem} from finitely many to countably many points would be to show that all of the intermediate results we use in this paper to prove Theorem~\ref{finitepointstheorem} hold uniformly in $n$.  As we prove these intermediate results, most of them will be easily seen to hold uniformly in $n$.  However, it is not clear whether or not part (iii) in Lemma~\ref{embeddingslemma} is uniform; it appears that the number $r_0$ which we obtain depends on $n$, and approaches 0 as $n$ tends to infinity.

We now return to the issue at hand: the finitely many points case.  Theorem~\ref{finitepointstheorem} says the following.  Given $n$
points on the $x_3-\mbox{axis}$, $(0,0,b_j)$ for $j=1,\ldots,n$,
with $b_1<\ldots <b_n$, we construct a sequence of disks
$\Sigma_i\subset B_1=B_1(0)\subset\mathbb{R}^3$ where the
curvatures blow up only at the prescribed $n$ points, and
$\Sigma_i\verb"\"\{x_3-\mbox{axis}\}$ consists of two multi-valued
graphs for each $i$.  The sequence
$\Sigma_i\verb"\"\cup_j\{x_3=b_j\}$ converges to $n+1$ embedded
minimal disks $\Sigma^k$, which sit between and spiral into the
appropriate planes $\{x_3=b_j\}$.  The result of Colding and
Minicozzi in \cite{CM1} is just Theorem~\ref{finitepointstheorem}
with $n=1$ and $b_1=0$.

For the reader's convenience, we will structure this paper
similarly to \cite{CM1}.  In particular, we provide some of the
brief background on the Weierstrass representation which Colding
and Minicozzi also outlined.

Let $\Omega\subset\mathbb{C}$ be a domain.  The Weierstrass
representation is as follows (see, for example, \cite{Os}).  Given
any meromorphic function $g$ on $\Omega$ and any holomorphic
one-form $\phi$ on $\Omega$, we obtain a (branched) conformal
minimal immersion $F:\Omega\rightarrow\mathbb{R}^3$, where
\begin{equation}\label{weierstrass}
    F(z)=\mbox{Re}\int_{\zeta\in\gamma_{z_0,z}}\left(\frac{1}{2}(g^{-1}(\zeta)-g(\zeta)),\frac{i}{2}(g^{-1}(\zeta)+g(\zeta)),1\right)\,\phi(\zeta).
\end{equation}
Here, we are integrating along a path $\gamma_{z_0,z}$ from a
fixed base point $z_0$ to $z$.  The choice of $z_0$ changes $F$ by
adding a constant.  We will assume that $F(z)$ is independent of
the choice of path, which is the case, for example, when $g$ has
no zeros or poles and $\Omega$ is simply connected (and this will
be the case for our choices of $g$ and $\Omega$).

The unit normal $\mathbf{n}$ and Gauss curvature $K$ of the
resulting minimal surface are given by (see \cite[Sec. 8,9]{Os})
\begin{equation}\label{unitnormal}
    \mathbf{n}=(2\,\mbox{Re}\,g,2\,\mbox{Im}\,g,|g|^2-1)/(|g|^2+1),
\end{equation}
\begin{equation}\label{gausscurv}
    K=-\left[\frac{4|\partial_z
    g||g|}{|\phi|(1+|g|^2)^2}\right]^2.
\end{equation}
The one-form $\phi$ is called the \textit{height differential},
and by equation~(\ref{unitnormal}), $g$ is the composition of the
Gauss map followed by stereographic projection.

We will assume that $\phi$ does not vanish and $g$ has no zeros or
poles; this implies that $F$ is an immersion, i.e., $dF\ne 0$. One
of the standard examples of this, which has the added benefit of
being an $\infty$-valued graph, and hence interesting for our
purposes, is the helicoid, whose Weierstrass data are
\begin{equation}\label{helicoiddata}
    g(z)=e^{iz}, \,\phi(z)=dz, \,\Omega=\mathbb{C}.
\end{equation}

This motivates the following.  If we want to construct
multi-valued minimal graphs, perhaps we should consider
Weierstrass data of the form
$$g(z)=e^{ih(z)}=e^{i(u(z)+iv(z))}, \,\phi(z)=dz,$$
for an appropriate choice of $\Omega$, where $h(z)$ is a
holomorphic function.  The next lemma gives us the differential of
$F$ in this case.

\begin{lemma}\label{differential}
If $F$ is given by equation~(\ref{weierstrass}) with
$g(z)=e^{i(u(z)+iv(z))}$ and $\phi(z)=dz$, then
\begin{equation}\label{partialx}
    \partial_x F=(\sinh v\,\cos u,\sinh v\,\sin u,1),
\end{equation}
\begin{equation}\label{partialy}
    \partial_y F=(\cosh v\,\sin u,-\cosh v\,\cos u,0).
\end{equation}
\end{lemma}

In particular, for the proof of Theorem~\ref{finitepointstheorem},
we will construct our multi-valued minimal graphs in this way,
with our choices of function $h_a(z)$ and domain $\Omega_a$
varying for each element of the sequence.  That is, we will
construct a one-parameter family of minimal immersions $F_a,
\,a\in (0,1/2)$, with Weierstrass data $g=e^{ih_a}$ (where
$h_a=u_a+iv_a$), $\phi =dz$, and domains $\Omega_a$ which we will
specify shortly.  We will prove that this family of immersions is
compact in Lemma~\ref{immersionscompactlemma}, and that the
immersions $F_a:\Omega_a\rightarrow\mathbb{R}^3$ are embeddings in
Lemma~\ref{embeddingslemma}.

For each $0<a<1/2$, let
\begin{equation}\label{defofhandomega}
    h_a(z)=\sum_{j=1}^n
    \frac{1}{2^{j-1}a}\,\mbox{arctan}\left(\frac{z-b_j}{a}\right)
    \,\,\mbox{on}\,\,
    \Omega_a=\cup_{j=1}^n\Omega_{a,j}, \,\,\mbox{where}
\end{equation}

\begin{eqnarray*}
    \Omega_{a,1} &=& \left\{(x,y):-\frac{1}{2}\le
    x\le\frac{b_2-b_1}{2},\,|y|\le\frac{[(x-b_1)^2+a^2]^{3/4}}{2}\right\}
    \nonumber\\
    \Omega_{a,j} &=& \left\{(x,y):\frac{b_j-b_{j-1}}{2}\le
    x\le\frac{b_{j+1}-b_j}{2},\,|y|\le\frac{[(x-b_j)^2+a^2]^{3/4}}{2}\right\},\nonumber\\
     & & {}j=2,\ldots,n-1 \nonumber\\
    \Omega_{a,n} &=& \left\{(x,y):\frac{b_n-b_{n-1}}{2}\le
    x\le\frac{1}{2},
    \,|y|\le\frac{[(x-b_n)^2+a^2]^{3/4}}{2}\right\}.
\end{eqnarray*}
To get an idea of what $\Omega_a$ looks like, note that the
$\Omega_{a,j}$ are defined similarly to the domain called
$\Omega_a$ by Colding and Minicozzi (see~\cite[Figure 4]{CM1}),
only centered at $b_j$ instead of at $0$.  When $a\rightarrow 0$,
the domain pinches off at the $n$ points $b_j$, just as Colding
and Minicozzi's domain pinches off at $0$ (see~\cite[Figure
5]{CM1}).

Note that $h_a$ is well-defined, since $\Omega_a$ is simply
connected and $b_j\pm ia\notin\Omega_a$ for $j=1,\ldots,n$.  By
direct computation, we see that
\begin{eqnarray}
    \partial_z h_a(z) &=&
    \sum_{j=1}^n\frac{1}{2^{j-1}}\,\frac{1}{(z-b_j)^2+a^2}\nonumber\\
      &=&
      \sum_{j=1}^n\frac{1}{2^{j-1}}\,\frac{(x-b_j)^2+a^2-y^2-2i(x-b_j)y}{[(x-b_j)^2+a^2-y^2]^2+4(x-b_j)^2y^2}.\label{partialz}
\end{eqnarray}
By the Cauchy-Riemann equations, we get
\begin{equation}\label{CRequations}
    \partial_z h_a = \partial_x u_a-i\partial_y u_a=\partial_y
    v_a+i\partial_x v_a.
\end{equation}
Also, the curvature is given by (see equation~(\ref{gausscurv}))
\begin{eqnarray}
    K_a(z) &=& \frac{-|\partial_z h_a|^2}{\cosh^4 v_a} \nonumber\\
      &=&-\frac{|\sum_{j=1}^n 2^{1-j}((z-b_j)^2+a^2)^{-1}|^2}{\cosh^4(\mbox{Im}(\sum_{j=1}^n\,\mbox{arctan}((z-b_j)/a)/2^{j-1}a))}.\label{gausscurvinourcase}
\end{eqnarray}
Note that $\lim_{a\rightarrow 0}|K_a(z)|=\infty$ for
$z=b_j,\,j=1,\ldots,n.$

Let $F_a:\Omega_a\rightarrow\mathbb{R}^3$ be from
equation~(\ref{weierstrass}) with $g=e^{ih_a}, \,\phi =dz$, and
$z_0=0$.  Let $\Omega_0=\cap_a\Omega_a\verb"\"\{b_1,\ldots,b_n\}$.
The family of functions $h_a$ is not compact, since
$\lim_{a\rightarrow 0}|h_a|(z)=\infty$ for $z\in\Omega_0$.
However, as the following lemma shows, the family of immersions
$F_a$ is compact.

\begin{lemma}\label{immersionscompactlemma}
If $a_k\rightarrow 0$, there exists a subsequence, which we also
call $a_k$, such that $F_{a_k}$ converges uniformly in $C^2$ on
compact subsets of $\Omega_0$.
\end{lemma}

\noindent{\it Proof}.  Similar to the proof of \cite[Lemma
2]{CM1}, with
$$-\sum_{j=1}^n\frac{1}{2^{j-1}}\,\frac{1}{z-b_j}$$
in place of -1/z.\qed

\bigskip In the next lemma, we show that the immersions
$F_a:\Omega_a\rightarrow\mathbb{R}^3$ are in fact embeddings. This
will follow from parts (i) and (ii) of the lemma.  Part (i) says
that the slice $\{x_3=t\}\cap F_a(\Omega_a)$ is the image of the
segment $\{x=t\}$ in the plane; that is, as $x$ varies and $y$
stays fixed, there is no self-intersection.  In part (ii), we show
that, in each slice $\{x_3=t\}\cap F_a(\Omega_a)$, the image
$F_a(\{x=t\}\cap\Omega_a)$ is a graph over some line segment in
the slice; that is, as $y$ varies and $x$ stays fixed, there is no
self-intersection.

\begin{lemma}\label{embeddingslemma}
For all $a>0$, the immersions
$F_a:\Omega_a\rightarrow\mathbb{R}^3$ satisfy
\begin{enumerate}
    \item[(i)] $x_3(F_a(x,y))=x$
    \item[(ii)] For each fixed $x$, $F_a(x,\cdot)$ is a graph in
    the plane $\{x_3=x\}$.
    \item[(iii)] There exists $r_0>0$ such that, for all $a$,
    \begin{eqnarray*}
    \left|F_a\left(x,\pm\frac{[(x-b_1)^2+a^2]^{3/4}}{2}\right)-F_a(x,0)\right|
    &>& r_0,\,-\frac{1}{2}\le x\le\frac{b_2-b_1}{2}\\
    \left|F_a\left(x,\pm\frac{[(x-b_j)^2+a^2]^{3/4}}{2}\right)-F_a(x,0)\right|
    &>& r_0,\,\frac{b_j-b_{j-1}}{2}\le
    x\le\frac{b_{j+1}-b_j}{2},\\
      & & {}j=2,\ldots,n-1\\
    \left|F_a\left(x,\pm\frac{[(x-b_n)^2+a^2]^{3/4}}{2}\right)-F_a(x,0)\right|
    &>& r_0,\,\frac{b_n-b_{n-1}}{2}\le x\le\frac{1}{2}.
    \end{eqnarray*}
\end{enumerate}
\end{lemma}

\medskip\noindent{\it Proof}.  (i) is immediate by the definition
of $F_a$, since $z_0=0$ and $\phi=dz$.

To prove (ii), first note that, by equations~(\ref{partialz}) and
(\ref{CRequations}), we have
$$|\partial_y
u_a(x,y)|\le\sum_{j=1}^n\frac{1}{2^{j-1}}\,\frac{2|x-b_j||y|}{[(x-b_j)^2+a^2-y^2]^2+4(x-b_j)^2y^2}.$$
Fix $k,\,1\le k\le n.$  On $\Omega_{a,k}$ (where
$(x-b_k)^2=\min_j(x-b_j)^2$), we have, for all $j=1,\ldots,n$,
\begin{eqnarray*}
    [(x-b_j)^2+a^2-y^2]^2+4(x-b_j)^2y^2 &\ge&
    [(x-b_j)^2+a^2-y^2]^2\\
     &\ge& \left[(x-b_j)^2+a^2-\frac{(x-b_k)^2+a^2}{4}\right]^2\\
     &\ge& \left[(x-b_j)^2+a^2-\frac{(x-b_j)^2+a^2}{4}\right]^2\\
     &=& \frac{9}{16}[(x-b_j)^2+a^2]^2.
\end{eqnarray*}
Therefore, we have
\begin{equation}\label{upartialyinequality}
    |\partial_y
    u_a(x,y)|\le 4\sum_{j=1}^n\frac{1}{2^{j-1}}\,\frac{|x-b_j||y|}{[(x-b_j)^2+a^2]^2}.
\end{equation}
Set $y_{x,a,k}=\frac{[(x-b_k)^2+a^2]^{3/4}}{2}$.
Integrating~(\ref{upartialyinequality}) gives
\begin{eqnarray}
    \max_{|y|\le y_{x,a,k}}|u_a(x,y)-u_a(x,0)| &\le& \max_{|y|\le
    y_{x,a,k}}\left|\int_0^y\partial_y u_a(x,t)\,dt\right|\nonumber\\
     &\le&
     \int_0^{y_{x,a,k}} 4\sum_{j=1}^n\frac{1}{2^{j-1}}\,\frac{|x-b_j|t}{[(x-b_j)^2+a^2]^2}\,dt\nonumber\\
     &=&
     2\sum_{j=1}^n\frac{1}{2^{j-1}}\,\frac{|x-b_j|t^2}{[(x-b_j)^2+a^2]^2}\left|_0^{y_{x,a,k}}\right.\nonumber\\
     &=&
     2\sum_{j=1}^n\frac{1}{2^{j-1}}\,\frac{|x-b_j|}{[(x-b_j)^2+a^2]^2}\,\frac{[(x-b_k)^2+a^2]^{3/2}}{4}\nonumber\\
     &\le&
     \frac{1}{2}\sum_{j=1}^n\frac{1}{2^{j-1}}\,\frac{|x-b_j|}{[(x-b_j)^2+a^2]^2}\,[(x-b_j)^2+a^2]^{3/2}\nonumber\\
     &=&
     \frac{1}{2}\sum_{j=1}^n\frac{1}{2^{j-1}}\,\frac{|x-b_j|}{[(x-b_j)^2+a^2]^{1/2}}\nonumber\\
     &\le& \frac{1}{2}\sum_{j=1}^n\frac{1}{2^{j-1}}\nonumber\\
     &<& \frac{1}{2}\sum_{j=1}^\infty \frac{1}{2^{j-1}}\nonumber\\
     &=& 1. \label{useparationinequality}
\end{eqnarray}
Set $\gamma_{x,a}(y)=F_a(x,y)$.  Since $v_a(x,0)=0$ and $\cos
(1)>1/2$, combining~(\ref{partialy}) and
(\ref{useparationinequality}), we get
\begin{eqnarray}
    \langle\gamma_{x,a}'(y),\gamma_{x,a}'(0)\rangle &=& \cosh
    v_a(x,y)\,\cos(u_a(x,y)-u_a(x,0))\nonumber\\
     &>& \cosh v_a(x,y)/2,\label{gammaprimeangle}
\end{eqnarray}
where $\gamma_{x,a}'(y)=\partial_y F_a(x,y)$.
By~(\ref{gammaprimeangle}), the angle between $\gamma_{x,a}'(y)$
and $\gamma_{x,a}'(0)$ is always less than $\pi /2$, proving (ii)
on $\Omega_{a,k}$, and hence on all of $\Omega_a$ since $k$ was
arbitrary.

To prove (iii), note that, by~(\ref{partialz}) and
(\ref{CRequations}), we have
$$\partial_y
v_a(x,y)=\sum_{j=1}^n\frac{1}{2^{j-1}}\,\frac{(x-b_j)^2+a^2-y^2}{[(x-b_j)^2+a^2-y^2]^2+4(x-b_j)^2y^2}.$$
As before, fix $k, \,1\le k\le n$, and look on $\Omega_{a,k}$.
Then, for all $j=1,\ldots,n$,
\begin{eqnarray*}
    [(x-b_j)^2+a^2-y^2]^2+4(x-b_j)^2y^2 &\le&
    [(x-b_j)^2+a^2+y^2]^2\\
      &\le& \left[(x-b_j)^2+a^2+\frac{(x-b_k)^2+a^2}{4}\right]^2\\
      &\le& \left[(x-b_j)^2+a^2+\frac{(x-b_j)^2+a^2}{4}\right]^2\\
      &=& \frac{25}{16}[(x-b_j)^2+a^2]^2.
\end{eqnarray*}
\begin{eqnarray*}
    (x-b_j)^2+a^2-y^2 &\ge&
    (x-b_j)^2+a^2-\frac{(x-b_k)^2+a^2}{4}\\
      &\ge& (x-b_j)^2+a^2-\frac{(x-b_j)^2+a^2}{4}\\
      &=& \frac{3}{4}[(x-b_j)^2+a^2].
\end{eqnarray*}
So, we have
\begin{eqnarray}
    \partial_y v_a(x,y) &\ge&
    \frac{12}{25}\sum_{j=1}^n\frac{1}{2^{j-1}}\,\frac{1}{(x-b_j)^2+a^2}\nonumber\\
      &>&
      \frac{3}{8}\sum_{j=1}^n\frac{1}{2^{j-1}}\,\frac{1}{(x-b_j)^2+a^2}.\label{vpartialyinequality}
\end{eqnarray}
Let $y_{x,a,k}=\frac{[(x-b_k)^2+a^2]^{3/4}}{2}$, as before.  Since
$v_a(x,0)=0$, integrating~(\ref{vpartialyinequality}) gives
\begin{eqnarray}
    \min_{\frac{y_{x,a,k}}{2}\le |y|\le y_{x,a,k}}|v_a(x,y)| &=&
    \min_{\frac{y_{x,a,k}}{2}\le |y|\le
    y_{x,a,k}}\left|\int_0^y\partial_y
    v_a(x,t)\,dt\right|\nonumber\\
      &>&
      \frac{3}{8}\int_0^{\frac{y_{x,a,k}}{2}}\sum_{j=1}^n\frac{1}{2^{j-1}}\,\frac{1}{(x-b_j)^2+a^2}\,dt\nonumber\\
      &\ge&
      \frac{3}{8}\int_0^{\frac{y_{x,a,k}}{2}}\frac{1}{2^{k-1}}\,\frac{1}{(x-b_k)^2+a^2}\,dt\nonumber\\
      &=&
      \frac{3}{8}\,\frac{1}{2^{k-1}}\,\frac{1}{(x-b_k)^2+a^2}\,\frac{[(x-b_k)^2+a^2]^{3/4}}{4}\nonumber\\
      &=& \frac{3}{32}\,\frac{1}{2^{k-1}}\,[(x-b_k)^2+a^2]^{-1/4}\nonumber\\
      &>&
      \frac{[(x-b_k)^2+a^2]^{-1/4}}{11\cdot 2^{n-1}}.\label{minvinequality}
\end{eqnarray}
Now, integrating~(\ref{gammaprimeangle}) and using
(\ref{minvinequality}), we obtain
\begin{equation}\label{gammaseparationinequality}
    \langle\gamma_{x,a}(y_{x,a,k})-\gamma_{x,a}(0),\gamma_{x,a}'(0)\rangle
    > \frac{[(x-b_k)^2+a^2]^{3/4}}{16}
    \,e^{[(x-b_k)^2+a^2]^{-1/4}/11\cdot 2^{n-1}}.
\end{equation}
Since $\lim_{s\rightarrow 0} s^3\,e^{s^{-1}/11\cdot
2^{n-1}}=\infty$, (\ref{gammaseparationinequality}) and its analog
for $\gamma_{x,a}(-y_{x,a,k})$ give an $r_k>0$ for which (iii)
holds (with $r_k$ in place of $r_0$) on $\Omega_{a,k}$.  This
proves (iii) on all of $\Omega_a$, with $r_0 =\min_k r_k$.\qed

\begin{corollary}\label{embeddingscorollary}
Let $r_0$ be given by part (iii) of Lemma~\ref{embeddingslemma}.
Then,
\begin{enumerate}
    \item[(a)] $F_a$ is an embedding.
    \item[(b)] $F_a(t,0)=(0,0,t)$ for $|t|<1/2$.
    \item[(c)] $\{0<x_1^2+x_2^2<r_0^2\}\cap
    F_a(\Omega_a)=\widetilde\Sigma_{1,a}\cup\widetilde\Sigma_{2,a}$
    for multi-valued graphs
    $\widetilde\Sigma_{1,a},\,\widetilde\Sigma_{2,a}$ over
    $D_{r_0}\verb"\"\{0\}$.
\end{enumerate}
\end{corollary}

\noindent {\it Proof}.  Same as \cite[Cor. 1]{CM1}.\qed

\bigskip\noindent {\it Proof of
Theorem~\ref{finitepointstheorem}}.  By scaling, it suffices to
find a sequence $\Sigma_i\subset B_R$ for some $R>0$.  By
Corollary~\ref{embeddingscorollary}, there exist minimal
embeddings $F_a:\Omega_a\rightarrow\mathbb{R}^3$ with
$F_a(t,0)=(0,0,t)$ for $|t|<1/2$, so (iii) holds for any $R\le
r_0$.  Set $R=\min\{r_0/2,1/4\}$, and $\Sigma_i=B_R\cap
F_{a_i}(\Omega_{a_i})$, where the sequence $a_i$ is to be
determined.

For each $j=1,\ldots,n$, by equation~(\ref{gausscurvinourcase}),
we have $|K_a|(b_j)\rightarrow\infty$ as $a\rightarrow 0$, proving
(i).

Also by~(\ref{gausscurvinourcase}), for each $j=1,\ldots,n$ and
all $\delta >0$,
$$\sup_a\sup_{\{|x-b_j|\ge\delta\}\cap\Omega_a} |K_a|<\infty$$
for all $x\notin\{b_1,...,b_n\}$.  Combined with (iii) and Heinz's
curvature estimate for minimal graphs (see, for example,
\cite[11.7]{Os}), this proves (ii).

By Lemma~\ref{immersionscompactlemma}, we can choose
$a_i\rightarrow 0$ so that the $F_{a_i}$ converge uniformly in
$C^2$ on compact subsets to $F_0:\Omega_0\rightarrow\mathbb{R}^3$.
So, by Lemma~\ref{embeddingslemma}, we obtain (iv)(a) and the
decomposition
$\Sigma^k\verb"\"\{x_3-\mbox{axis}\}=\Sigma_1^k\cup\Sigma_2^k$ for
multi-valued graphs $\Sigma_j^k$, where $j=1,2$ and
$k=1,\ldots,n+1$.  To obtain (iv)(b) and the remainder of (iv)(c),
we must show that each graph $\Sigma_j^k$ is $\infty$-valued, as
this would imply the spiraling which we seek.  By (iii) and
(\ref{partialy}), the level sets $\{x_3=x\}\cap\Sigma_j^k$ are
graphs over the line in the direction
$$\lim_{a\rightarrow 0}(\sin u_a(x,0),-\cos u_a(x,0),0).$$
Since, for all $j=1,\ldots,n$ and all $t$ sufficiently close to
$b_j$,
$$\lim_{a\rightarrow
0}|u_a(t-b_j,0)-u_a(2(t-b_j),0)|=\frac{1}{2(t-b_j)},$$ we see
that, for $t$ sufficiently close to $b_j$,
$\{t-b_j<|x_3|<2(t-b_j)\}\cap\Sigma_j^k$ contains an embedded
$N_t$-valued graph, where $N_t\approx 1/4\pi
(t-b_j)\rightarrow\infty$ as $t\rightarrow b_j$. This proves that
each $\Sigma_j^k$ spirals the way we claim, completing the proof
of (iv).\qed

\end{document}